\documentclass[12pt]{article}
\usepackage{amssymb}
\usepackage{amsfonts}
\usepackage{amsmath}
\usepackage{amsbsy}
\newcommand{\be}{\begin{equation}}
\newcommand{\ee}{\end{equation}}
\newcommand{\bea}{\begin{eqnarray}}
\newcommand{\eea}{\end{eqnarray}}
\newcommand{\ba}{\begin{array}}
\newcommand{\ea}{\end{array}}

\newcommand{\bc}{\begin{center}}
\newcommand{\ec}{\end{center}}
\newcommand{\ben}{\begin{enumerate}}
\newcommand{\een}{\end{enumerate}}
\newcommand{\bfi}{\begin{figure}}
\newcommand{\efi}{\end{figure}}

\newcommand{\bq}{\begin{quote}}
\newcommand{\eq}{\end{quote}}
\newcommand{\bqu}{\begin{quotation}}
\newcommand{\equ}{\end{quotation}}
\newenvironment{emphit}{\begin{itemize}}{\end{itemize}}
\newcommand{\bemp}{\begin{emphit}}
\newcommand{\eemp}{\end{emphit}}

\newcommand{\bt}{\begin{tabular}}
\newcommand{\et}{\end{tabular}}

\newtheorem{myth}{Theorem}[section]
\newtheorem{mylem}{Lemma}[section]
\newtheorem{mycor}{Corollary}[section]

\newtheorem{mydef}{Definition}[section]

\begin{document}
\date{}
\title{On the existence of A-loops with some commutative inner mappings and others of order
$2$
}
\author{T. G. Jaiy\'e\d ol\'a\\
Department of Mathematics,\\
Obafemi Awolowo University,\\
Ile Ife 220005, Nigeria.\\
jaiyeolatemitope@yahoo.com\\tjayeola@oauife.edu.ng \and
J. O. Ad\'en\'iran\thanks{All correspondence to be addressed to this author.}  \\
Department of Mathematics,\\
University of Abeokuta, \\
Abeokuta 110101, Nigeria.\\
ekenedilichineke@yahoo.com\\
adeniranoj@unaab.edu.ng} \maketitle

\begin{abstract} The existence
of A$_\rho$-loops, A$_\lambda$-loops and A$_\mu$-loops that are
neither extra loops nor CC-loops such that any two of their inner
mappings $R(x,y),L(x,y)$ and $T(x)$ commute while the other one is
of order $2$ is shown.
\end{abstract}

\section{\textbf{INTRODUCTION}}
Let $L$ be a non-empty set. Define a binary operation ($\cdot $) on
$L$. If $x\cdot y\in L$ for all $x, y\in L$, $(L, \cdot )$ is called
a groupoid. If the system of equations ; $a\cdot x=b$ and $y\cdot
a=b$ have unique solutions for $x$ and $y$ respectively, then $(L,
\cdot )$ is called a quasigroup. Furthermore, if there exists a
unique element $e\in L$ called the identity element such that for
all $x\in L$, $x\cdot e=e\cdot x=x$, $(L, \cdot )$ is called a loop.
A detailed information on loop properties, types, concepts and
applications are contained in \cite{phd3}, \cite{phd41},
\cite{phd39}, \cite{phd49}, \cite{phd42} and \cite{phd75}.

The symmetric group of a loop $(L, \cdot )$ is denoted by $S(L,
\cdot )$ and it is defined as the group of all permutations or
self-bijections on $L$. The bijection $L_x : L\to L$ defined as
$yL_x=x\cdot y$ for all $x, y\in L$ is called a left
translation(multiplication) of $L$ while the bijection $R_x : L\to
L$ defined as $yR_x=y\cdot x$ for all $x, y\in L$ is called a right
translation(multiplication) of $L$. In a loop $(L,\cdot )$, the
group generated by the set of left or right translations and their
inverses is denoted by ${\cal M}_\lambda (L,\cdot )$ or ${\cal
M}_\rho (L,\cdot )$ and called the left or right multiplication
group of $(L,\cdot )$, while the group generated by the set of both
left and right translations and their inverses is denoted by ${\cal
M} (L,\cdot )$ and called the multiplication group of $(L,\cdot )$.
It is well known that the groups ${\cal M}_\lambda (L,\cdot ),{\cal
M}_\rho (L,\cdot )$ and ${\cal M}(L,\cdot )$ are subgroups of
$S(L,\cdot )$.

The triple $(U, V, W)$ formed such that $U, V, W\in S(L, \cdot )$ is
called an autotopism of $L$ if and only if $xU\cdot yV=(x\cdot y)W$
for all $x,y\in L$. The group of autotopisms of $L$ is called the
autotopism group and it is denoted by $AUT(L, \cdot )$. If $U=V=W$,
then $U$ is called an automorphism. The group of automorphisms on
$L$ is called the automorphism group and it is denoted by $A(L,
\cdot )$. If $U\in S(L, \cdot )$ such that $(U,UR_c,UR_c)\in AUT(L,
\cdot)$ or $(UL_c,U,UL_c)\in AUT(L, \cdot )$ for some $c\in L$, then
$U$ is called a right or left pseudo-automorphism of the loop $L$
with a right or left companion $c$. The group formed by such
permutations is called the right or left pseudo-automorphism group
and it is denoted by $PS_\rho(L,\cdot )$ or $PS_\lambda(L,\cdot )$.

All elements $\alpha \in {\cal M}_\lambda (L,\cdot )$ or $\alpha \in
{\cal M}_\rho (L,\cdot )$ or $\alpha \in {\cal M}(L,\cdot )$  such
that $e\alpha =e$ form a group called the left inner mapping group
or right inner mapping group or inner mapping group of $(L,\cdot )$
and this is denoted by $\textrm{Inn}_\lambda(L,\cdot )$ or
$\textrm{Inn}_\rho (L,\cdot )$ or $\textrm{Inn}(L,\cdot )$.

The inner mapping $R(x,y)=R_xR_yR_{xy}^{-1}$ is called a right inner
mapping and it has been shown that they generate the group
$\textrm{Inn}_\rho (L,\cdot )$. A loop $L$ is called a right
A-loop(A$_\rho$-loop) if $\textrm{Inn}_\rho (L,\cdot )\le A(L,\cdot
)$.

The inner mapping $L(x,y)=L_xL_yL_{yx}^{-1}$ is called a left inner
mapping and it has been shown that they generate the group
$\textrm{Inn}_\lambda(L,\cdot )$. A loop $L$ is called a left
A-loop(A$_\lambda$-loop) if $\textrm{Inn}_\lambda(L,\cdot )\le
A(L,\cdot )$.

The inner mapping $T(x)=R_xL_x^{-1}$ is called a middle inner
mapping and the group generated by these is denoted by
$\textrm{Inn}_\mu (L,\cdot )$ and called the middle inner mapping
group. A loop $L$ is called a middle A-loop(A$_\mu$-loop) if
$\textrm{Inn}_\mu (L,\cdot )\le A(L,\cdot )$.

It has been shown in \cite{phd3} that the inner mapping group
$\textrm{Inn}(L,\cdot )$ of a loop $L$ is generated by its left,
right and middle inner mappings. So, if $\textrm{Inn}(L,\cdot )\le
A(L,\cdot )$, $L$ is called an A-loop, hence, $L$ is an A-loop if
and only if $L$ is an A$_\rho$-loop, A$_\lambda$-loop and an
A$_\mu$-loop. The study of A-loops started by Bruck and Paige in
\cite{phd40}. Further studies on A-loops have been done by Osborn
\cite{phd110}, Phillips \cite{phd123} and Drapal \cite{phd14}. The
most interesting work on A-loops is Kinyon et. al. \cite{phd13}
which gives the solution to the Osborn problem.

After the introduction of conjugacy closed loops(CC-loop) by
Goodaire and Robinson \cite{phd91}, \cite{phd48}, a tremendous study
of their properties and structural behaviours have been studied by
Kunen \cite{phd78} and some recent works of Kinyon and Kunen
\cite{phd36}, \cite{phd47}, Phillips et. al. \cite{phd35}, Dr\'apal
\cite{phd37}, \cite{phd38}, \cite{phd98}, \cite{phd139},
\cite{phd107}, Cs\"org\H o et. al. \cite{phd106}, \cite{phd104},
\cite{phd108} and Phillips \cite{phd151}. In \cite{phd36},
\cite{phd35}, \cite{phd47} and \cite{phd78}, it is proved and stated
that in a CC-loop $L$;
\begin{itemize}
\item $R(x,y),L(x,y)\in A(L)$, hence $L$ is a both an A$_\rho$-loop and an
A$_\lambda$-loop,
\item $\textrm{Inn}_\lambda(L)=\textrm{Inn}_\rho (L)$ and $\textrm{Inn}(L,\cdot )=\Big<\{T(x)~:~x\in L\}\Big>$,
\item $R(x,y)R(u,v)=R(u,v)R(x,y)$ and $R(x,y)L(u,v)=L(u,v)R(x,y)$,
hence $\textrm{Inn}_\rho (L)$ and $\Big<\{R(x,y),L(x,y)~:~x,y\in
L\}\Big>$ are abelian groups.
\end{itemize}
If $L$ is an extra loop, then the facts listed above and the ones
below are true.
\begin{itemize}
\item $R(x,y)=L(x,y)=R(y,x)=L(y,x)$, $|R(x,y)|=2$, hence $\textrm{Inn}_\lambda(L)=\textrm{Inn}_\rho
(L)$ are boolean groups,
\item $T(x)\in A(L)$ if and only if $x\in N(L)$.
\end{itemize}
If $L$ is an A-loop then according to \cite{phd40},
$T(x)L(y,x)=L(y,x)T(x)$ and $T(x)R(x,y)=R(x,y)T(x)$.

The multiplication group and inner mapping group of loops have been
studied by Dr\'apal \cite{phd100}, \cite{phd101}, \cite{phd102},
\cite{phd99}, Dr\'apal et. al. \cite{phd129}, \cite{phd130}, Kepka
\cite{phd131}, \cite{phd132}, Kepka and Niemenmaa \cite{phd133},
Niemenmaa \cite{phd135}, \cite{phd136}, \cite{phd162}, Niemenmaa and
Kepka \cite{phd134}, Cs\"org\H o and Kepka \cite{phd103} in
different fashions. The multiplication group structure determines
the structure of a loop(e.g solvability of ${\cal M}(L)$ implies the
solvability of a finite loop $L$ \cite{phd163}) while if
$\textrm{Inn}(L)$ is of order $2p$($p$ an odd prime), then ${\cal
M}(L)$ is solvable hence $L$ is solvable as well(\cite{phd164}).
\paragraph{}
The present study investigates the existence of A$_\rho$-loops,
A$_\lambda$-loops, A$_\mu$-loops and A-loops that are neither extra
loops nor CC-loops such that any two of their inner mappings
$R(x,y),L(x,y)$ and $T(x)$ commute while the other one is of order
$2$.

\begin{mydef}
If $(L, \cdot )$ and $(G, \circ )$ are two distinct loops, then the
triple $(U, V, W) : (L, \cdot )\rightarrow (G, \circ )$ such that
$U, V, W : L\rightarrow G$ are bijections is called a loop isotopism
if and only if
\begin{displaymath}
xU\circ yV=(x\cdot y)W~\forall ~x, y\in L.
\end{displaymath}
Throughout, when $L_x~:~y\mapsto xy$ and $R_x~:~y\mapsto yx$ are
respectively the left and right translations of a loop then the left
and right translations of its loop isotope are denoted by
$L_x'~:~y\mapsto xy$ and $R_x'~:~y\mapsto yx$ respectively.
\end{mydef}

\begin{mydef}\label{definition:deviation}
Let $(B,\cdot )$ be a loop. If $x\in B$ and $\phi\in S(B,\cdot )$,
then the mapping
\begin{displaymath}
\mu_x(\phi)~:~S(B,\cdot)\longrightarrow S(B,\cdot)~\textrm{defined
by}~\mu_x(\phi)=\phi^{-1} L_x\phi L_{x\phi }^{-1}
\end{displaymath}
is called the deviation of the mapping $\phi$ at $x$.\\

Furthermore, set
\begin{displaymath}
P(x,\phi ):=L_x\phi -\phi L_{x\phi }^{-1}\phi L_x\phi^{-1}L_{x\phi
}.
\end{displaymath}
\end{mydef}

\section{\textbf{MAIN RESULTS}}
\subsection{Deviation}
\begin{mylem}\label{1:1}
Let $(B,\cdot )$ be a loop with $\phi\in S(B,\cdot )$.
\begin{displaymath}
\textrm{If}~P(x,\phi ):=0,~\textrm{then}~\mu_x(\phi)=L_{x\phi
}^{-1}\phi L_x\phi^{-1}.
\end{displaymath}
\end{mylem}
\textit{Proof} \\
If $P(x,\phi ):=0$, then $L_x\phi =\phi L_{x\phi }^{-1}\phi
L_x\phi^{-1}L_{x\phi }\Longrightarrow \phi ^{-1}L_x\phi L_{x\phi
}^{-1} =L_{x\phi }^{-1}\phi L_x\phi^{-1}\Longrightarrow
\mu_x(\phi)=L_{x\phi }^{-1}\phi L_x\phi^{-1}$ since
$\mu_x(\phi)=\phi^{-1} L_x\phi L_{x\phi }^{-1}$.

\begin{myth}\label{1:2}
Let $(B,\cdot )$ be a loop with $\phi\in S(B,\cdot )$ such that
$\phi~:~e\mapsto e$. If $P(x,\phi ):=0$, then $\phi\in A(B,\cdot )$.
\end{myth}
\textit{Proof}\\
$P(x,\phi ):=0\Longrightarrow L_x\phi =\phi L_{x\phi }^{-1}\phi
L_x\phi^{-1}L_{x\phi }\Longrightarrow yL_x\phi =y\phi L_{x\phi
}^{-1}\phi L_x\phi^{-1}L_{x\phi }\Longrightarrow $
\begin{equation}\label{eq:1}
(xy)\phi =y\phi L_{x\phi }^{-1}\phi L_x\phi^{-1}L_{x\phi }.
\end{equation}
Let $z=y\phi L_{x\phi }^{-1}\Longrightarrow x\phi\cdot z=y\phi$. Put
$x=e$, then
\begin{equation}\label{eq:2}
z=y\phi.
\end{equation}
Now from equation (\ref{eq:1}), we have $(xy)\phi =z\phi
L_x\phi^{-1}L_{x\phi }=(x\cdot z\phi )\phi^{-1}L_{x\phi }=x\phi\cdot
(x\cdot z\phi )\phi^{-1}$. So,
\begin{equation}\label{eq:3}
(xy)\phi =x\phi\cdot (x\cdot z\phi )\phi^{-1}.
\end{equation}
Let $z'=(x\cdot z\phi )\phi^{-1}$, then $x\cdot z\phi =z'\phi$.
Using equation (\ref{eq:2}), $x\cdot y\phi^2=z'\phi$. Now, let
$x=e$, then $y\phi^2=z'\phi\Longrightarrow $
\begin{equation}\label{eq:4}
z'=y\phi.
\end{equation}
From equation (\ref{eq:3}), $(xy)\phi =x\phi\cdot z'$. So by
equation (\ref{eq:4}),
\begin{displaymath}
(xy)\phi =x\phi\cdot y\phi\Longrightarrow \phi\in A(B,\cdot ).
\end{displaymath}

\begin{myth}\label{automorphism:exponent2}
Let $(B,\cdot )$ be a loop and $\phi\in S(B,\cdot )$. The following
are true.
\begin{enumerate}
\item $\phi~:~e\mapsto e\Leftrightarrow \mu_x(\phi)~:~e\mapsto e~\forall~x\in B$.
\item $\phi\in PS_\lambda(B,\cdot )\Leftrightarrow ~\exists~c\in B~\ni~\mu_x(\phi)=L(x\phi,e)^{-1}~\forall~x\in B$.
\item $\phi\in A(B,\cdot )\Leftrightarrow \mu_x(\phi)=I~\forall~x\in
B$.
\end{enumerate}
\end{myth}
\textit{Proof}\\
Recall that $\mu_x(\phi)=\phi^{-1} L_x\phi L_{x\phi }^{-1}$.
\begin{enumerate}
\item $e\mu_x(\phi)=e\Longleftrightarrow e\phi^{-1} L_x\phi L_{x\phi
}^{-1}=e\Longleftrightarrow e\phi^{-1}L_x\phi =eL_{x\phi
}\Longleftrightarrow (x\cdot e\phi^{-1})\phi =x\phi
\Longleftrightarrow x\cdot e\phi^{-1}=x\Longleftrightarrow
e\phi^{-1}=e\Longleftrightarrow e\phi=e$.
\item $\phi\in PS_\lambda(B,\cdot )$ with a left companion $c\in
B$ if and only if
\begin{equation}\label{eq:5}
c\cdot (x\cdot y)\phi =(c\cdot x\phi )\cdot y\phi\Longleftrightarrow
L_{(c\cdot x\phi )}=\phi^{-1}L_x\phi L_c.
\end{equation}
$L(x,y)=L_xL_yL_{yx}^{-1}$, so $L(x\phi,c)=L_{x\phi}L_cL_{(c\cdot
x\phi)}^{-1}$. Thus, computing and using equation (\ref{eq:5}),
\begin{displaymath}
\mu_x(\phi)L(x\phi,c)=\Big(\phi^{-1} L_x\phi L_{x\phi
}^{-1}\Big)\Big(L_{x\phi}L_cL_{(c\cdot x\phi)}^{-1}\Big)=\phi^{-1}
L_x\phi L_cL_{(c\cdot x\phi)}^{-1}=L_{(c\cdot x\phi)}L_{(c\cdot
x\phi)}^{-1}=I
\end{displaymath}
which implies $\mu_x(\phi)L(x\phi,c)=I\Longrightarrow
\mu_x(\phi)=L(x\phi,c)^{-1}$.\\

Conversely, if $\mu_x(\phi)=L(x\phi,c)^{-1}$, then
$\mu_x(\phi)L(x\phi,c)=I$. So,
\begin{displaymath}
\Big(\phi^{-1} L_x\phi L_{x\phi
}^{-1}\Big)\Big(L_{x\phi}L_cL_{(c\cdot
x\phi)}^{-1}\Big)=I\Longrightarrow \phi^{-1} L_x\phi L_cL_{(c\cdot
x\phi)}^{-1}=I\Longrightarrow \phi^{-1} L_x\phi L_c=L_{(c\cdot
x\phi)}
\end{displaymath}
implies $\phi\in PS_\lambda(B,\cdot )$ with a left companion $c\in
B$ by following equation (\ref{eq:5}).
\item Following 2., $\phi\in A(B,\cdot )\Longleftrightarrow \phi\in PS_\lambda(B,\cdot )$ such that
$c=e$. $L(x\phi,e)=I$,
\begin{displaymath}
\therefore~\phi\in A(B,\cdot )\Leftrightarrow
\mu_x(\phi)=I~\forall~x\in B.
\end{displaymath}
\end{enumerate}

\begin{myth}\label{aut:exp2}
Let $(B,\cdot )$ be a loop with $\phi\in S(B,\cdot )$ such that
$\phi~:~e\mapsto e$. If $P(x,\phi ):=0$, then $|\phi|=2$.
\end{myth}
\textit{Proof}\\
By Theorem~\ref{1:2}, $\phi\in A(B,\cdot )$. So following
Theorem~\ref{automorphism:exponent2} and Lemma~\ref{1:1},
$\mu_x(\phi)=L_{x\phi }^{-1}\phi L_x\phi^{-1}=I$. Thus, $L_{x\phi
}^{-1}\phi L_x\phi^{-1}=I\Longrightarrow \phi L_x=L_{x\phi
}\phi\Longrightarrow y\phi L_x=yL_{x\phi }\phi\Longrightarrow x\cdot
y\phi=(x\phi\cdot y)\phi=x\phi^2\cdot y\phi\Longrightarrow x
=x\phi^2\Longrightarrow \phi^2=I$.

\subsection{Isotopic Characterization Of A-loops}
\begin{myth}\label{rita:iso}
Let $(G,\cdot )$ and $(H,\circ )$ be any two distinct quasigroups.
If $A,B,C : G\rightarrow H$ are permutations, then the following
conditions are equivalent :
\begin{enumerate}
\item the triple $\alpha=(A,B,C)$ is an isotopism of $G$ upon $H$.
\item $R_{xB}'=A^{-1}R_xC~\forall~x\in G$.
\item $L_{yA}'=B^{-1}L_yC~\forall~y\in G$.
\end{enumerate}
\end{myth}
\textit{Proof}\\
\begin{description}
\item[(1$\Leftrightarrow $2)] If $\alpha=(A,B,C) : (G,\cdot )\rightarrow (H,\circ )$ is an isotopism, then
$xA\circ yB=(x\cdot y)C\Leftrightarrow
xAR_{yB}'=xR_yC\Leftrightarrow AR_{yB}'=R_yC\Leftrightarrow
R_{yB}'=A^{-1}R_yC$.
\item[(1$\Leftrightarrow $3)] If
$\alpha=(A,B,C) : (G,\cdot )\rightarrow (H,\circ )$ is an isotopism,
then $xA\circ yB=(x\cdot y)C\Leftrightarrow
yBL_{xA}'=yL_xC\Leftrightarrow BL_{xA}'=L_xC\Leftrightarrow
L_{xA}'=B^{-1}L_xC$.
\end{description}
Finally, 1 $\Leftrightarrow $ 2 and 1 $\Leftrightarrow $ 3
$\Rightarrow $ 2 $\Leftrightarrow $ 3. Hence the statements 1, 2 and
3 are equivalent to each other.

\begin{myth}\label{main:result}
Let $G=(\Omega ,\cdot )$ and $G'=(\Omega ,\circ )$ be any two
distinct loops. If $A,B,C\in S(\Omega )$ such that
$P(x,\phi):=0~\forall~\phi\in \{A,B,C\}$, then the following are
equivalent for all $x\in \Omega $.
\begin{enumerate}
\item $(A,B,C)$ is an isotopism of $G$ upon $G'$.
\item $\mu_x(A)=CL_{xA^2}'^{-1}B^{-1}AL_xA^{-1}$.
\item $\mu_x(A)=L_{xA}^{-1}ABL_{xA}'(AC)^{-1}$.
\item $\mu_x(B)=CL_{xBA}'^{-1}L_xB^{-1}$.
\item $\mu_x(B)=L_{xB}^{-1}B^2L_{xA}'(BC)^{-1}$.
\item $\mu_x(C)=CL_{xCA}'^{-1}B^{-1}CL_xC^{-1}$.
\item $\mu_x(C)=L_{xC}^{-1}CBL_{xA}'C^{-2}$.
\item $\Big(I,AB,\mu_{xA^{-1}}(A)AC\Big)$ is an isotopism of $G$ upon $G'$.
\item $\Big(B^{-1}A,B^2,\mu_{xB^{-1}}(B)BC\Big)$ is an isotopism of $G$ upon $G'$.
\item $\Big(C^{-1}A,CB,\mu_{xC^{-1}}(C)C^2\Big)$ is an isotopism of $G$ upon $G'$.
\end{enumerate}
\end{myth}
\textit{Proof} \\
1 $\Leftrightarrow $ 2, 1 $\Leftrightarrow $ 3, 1 $\Leftrightarrow $
4, 1 $\Leftrightarrow $ 5, 1 $\Leftrightarrow $ 6, 1
$\Leftrightarrow $ 7, 1 $\Leftrightarrow $ 8, 1 $\Leftrightarrow $ 9
and 1 $\Leftrightarrow $ 10 are achieved by using
Theorem~\ref{rita:iso}, Lemma~\ref{1:1} and
Definition~\ref{definition:deviation}.

\begin{myth}\label{main:result1}
Let $G=(\Omega ,\cdot )$ and $G'=(\Omega ,\circ )$ be two distinct
loops.
\begin{enumerate}
\item If $A\in S(\Omega )$ such that $A~:~e\mapsto e$ and $P(x,A):=0~\forall~x\in \Omega$, then the following are equivalent:
\begin{description}
\item[(i)] $(A,B,C)$ is an isotopism of $G$ upon $G'$.
\item[(ii)] $(I,AB,AC)$ is an isotopism of $G$ upon $G'$.
\item[(iii)] $AL_{xA}AC=BL_{xA}'$.
\item[(iv)] $A=CL_x'B^{-1}AL_x$.
\end{description}
\item If $B\in S(\Omega )$ such that $B~:~e\mapsto e$ and $P(x,B):=0~\forall~x\in \Omega$, then the following are equivalent:
\begin{description}
\item[(i)] $(A,B,C)$ is an isotopism of $G$ upon $G'$.
\item[(ii)] $(BA,I,BC)$ is an isotopism of $G$ upon $G'$.
\item[(iii)] $B=CL_{xBA}'^{-1}L_x$.
\end{description}
\item If $C\in S(\Omega )$ such that $C~:~e\mapsto e$ and $P(x,C):=0~\forall~x\in \Omega$, then the following are equivalent:
\begin{description}
\item[(i)] $(A,B,C)$ is an isotopism of $G$ upon $G'$.
\item[(ii)] $(CA,CB,I)$ is an isotopism of $G$ upon $G'$.
\item[(iii)] $L_x=CBL_{xCA}'$.
\end{description}
\end{enumerate}
\end{myth}
\textit{Proof} \\
The proof lies wholly on Theorem~\ref{main:result}. And the outcomes
are achieved by using Theorem~\ref{automorphism:exponent2},
Theorem~\ref{aut:exp2} and Theorem~\ref{1:2}.

\begin{mycor}\label{main:result2}
Let $G=(\Omega ,\cdot )$ and $G'=(\Omega ,\circ )$ be two distinct
isotopic loops with different identity elements such that the triple
$(A,B,C)$ is the isotopism between $G$ and $G'$.
\begin{enumerate}
\item If $A\in S(\Omega )$ such that $A~:~e\mapsto e$ and $P(x,A):=0~\forall~x\in \Omega$,
$(A,B,C)$ is an isotopism of $G$ upon $G'$ if and only if
$L_x'=B^{-1}L_{xA}C$. Hence,
\begin{description}
\item[(i)] $C=BL_e'$, $B=CL_e'$ and $L_e'^2=I$.
\item[(ii)] $C^{-1}B=B^{-1}C$ and $CB=BC$.
\end{description}
\item If $B\in S(\Omega )$ such that $B~:~e\mapsto e$ and $P(x,B):=0~\forall~x\in
\Omega$, then $C=BL_{eA}'$.
\item If $C\in S(\Omega )$ such that $C~:~e\mapsto e$ and $P(x,C):=0~\forall~x\in \Omega$,
then $C=BL_{eA}'$.
\end{enumerate}
\end{mycor}
\textit{Proof} \\ The proof of this is a consequence of
Theorem~\ref{main:result1} by replacing $x\in \Omega$ with the
identity element $e$ of $G$.

\paragraph{}
The tables below summarize the important results of this subsection
as shown in Theorem~\ref{main:result1} and
Corollary~\ref{main:result2}.

\begin{center}
\begin{tabular}{|c|c|c|}
\hline
Hypothesis & Hypothesis & Inference \\
\hline
$P(x,\cdot ):=0$ & $\cdot~:~e\mapsto e$ & \\
\hline
A & A & $C=BL_e',B=CL_e'$  \\
  &   & $C^{-1}B=B^{-1}C,CB=BC$\\
\hline
B & B & $C=BL_e'$ \\
\hline
C & C & $C=BL_e'$ \\
\hline
\end{tabular}
\end{center}

\begin{center}
\begin{tabular}{|c|c|c|c|c|}
\hline
Hypothesis & Hypothesis & Hypothesis & Inference  \\
\hline
$\cdot~:~e\mapsto e$ & $P(x,\cdot ):=0$ & Isotopism & Equivalent Isotopism \\
\hline
A & A & (A,B,C) & (I,AB,AC) \\
\hline
B & B & (A,B,C) & (BA,I,BC) \\
\hline
C & C & (A,B,C) & (CA,CB,I) \\
\hline
\end{tabular}
\end{center}

\begin{myth}\label{2:1}
Let $G=(\Omega ,\cdot )$ and $G'=(\Omega ,\circ )$ be two distinct
isotopic loops.
\begin{enumerate}
\item Under the triple $\big(R(x,y),L(u,v),T(z)\big)$,
\begin{description}
\item[(a)] if $P(z,R(x,y)):=0$ then,
\begin{description}
\item[(i)] $G$ is an A$_\rho$-loop and $|R(x,y)|=2$.
\item[(ii)] $\textrm{Inn}_\mu (G)=\big<L(x,y)L_e'~:~x,y\in \Omega\big>$ and $\textrm{Inn}_\lambda (G)=\big<T(x)L_e'~:~x\in
\Omega\big>$.
\item[(iii)] $T(z)L(x,y)=L(x,y)T(z)$ and $T(z)^{-1}L(x,y)=L(x,y)^{-1}T(z)$, hence $L(x,y)^2=T(z)^2$.
\item[(iv)] the triple $\big(I,R(x,y)L(u,v),R(x,y)T(z)\big)$ is an isotopism
from $G$ to $G'$.
\end{description}
\item[(b)] if $P(z,L(x,y)):=0$ then,
\begin{description}
\item[(i)] $G$ is an A$_\lambda$-loop and $|L(x,y)|=2$.
\item[(ii)] $\textrm{Inn}_\mu (G)=\big<L(x,y)L_e'~:~x,y\in \Omega\big>$.
\item[(iii)] the triple $\big(L(u,v)R(x,y),I,L(u,v)T(z)\big)$ is an isotopism
from $G$ to $G'$.
\end{description}
\item[(c)] if $P(z,T(x)):=0$ then,
\begin{description}
\item[(i)] $G$ is an A$_\mu$-loop and $|T(x)|=2$.
\item[(ii)] $\textrm{Inn}_\mu (G)=\big<L(x,y)L_e'~:~x,y\in \Omega\big>$.
\item[(iii)] the triple $\big(T(z)R(x,y),T(z)L(u,v),I\big)$ is an isotopism
from $G$ to $G'$.
\end{description}
\end{description}
\item Under the triple $\big(L(x,y),R(u,v),T(z)\big)$,
\begin{description}
\item[(a)] if $P(z,L(x,y)):=0$ then,
\begin{description}
\item[(i)] $G$ is an A$_\lambda$-loop and $|L(x,y)|=2$.
\item[(ii)] $\textrm{Inn}_\mu (G)=\big<R(x,y)L_e'~:~x,y\in \Omega\big>$ and $\textrm{Inn}_\rho (G)=\big<T(x)L_e'~:~x\in
\Omega\big>$.
\item[(iii)] $T(z)R(x,y)=R(x,y)T(z)$ and $T(z)^{-1}R(x,y)=R(x,y)^{-1}T(z)$, hence $R(x,y)^2=T(z)^2$.
\item[(iv)] the triple $\big(I,L(x,y)R(u,v),L(x,y)T(z)\big)$ is an isotopism
from $G$ to $G'$.
\end{description}
\item[(b)] if $P(z,R(x,y)):=0$ then,
\begin{description}
\item[(i)] $G$ is an A$_\rho$-loop and $|R(x,y)|=2$.
\item[(ii)] $\textrm{Inn}_\mu (G)=\big<R(x,y)L_e'~:~x,y\in \Omega\big>$.
\item[(iii)] the triple $\big(R(x,y)L(u,v),I,R(x,y)T(z)\big)$ is an isotopism
from $G$ to $G'$.
\end{description}
\item[(c)] if $P(z,T(x)):=0$ then,
\begin{description}
\item[(i)] $G$ is an A$_\mu$-loop and $|T(x)|=2$.
\item[(ii)] $\textrm{Inn}_\mu (G)=\big<R(x,y)L_e'~:~x,y\in \Omega\big>$.
\item[(iii)] the triple $\big(T(z)L(x,y),T(z)R(u,v),I\big)$ is an isotopism
from $G$ to $G'$.
\end{description}
\end{description}
\item Under the triple $\big(T(z),R(x,y),L(u,v)\big)$,
\begin{description}
\item[(a)] if $P(y,T(x)):=0$ then,
\begin{description}
\item[(i)] $G$ is an A$_\mu$-loop and $|L(x,y)|=2$.
\item[(ii)] $\textrm{Inn}_\lambda (G)=\big<R(x,y)L_e'~:~x,y\in \Omega\big>$ and $\textrm{Inn}_\rho (G)=\big<L(x,y)L_e'~:~x\in
\Omega\big>$.
\item[(iii)] $R(x,y)L(u,v)=L(u,v)R(x,y)$ and $L(u,v)^{-1}R(x,y)=R(x,y)^{-1}L(u,v)$, hence $R(x,y)^2=L(u,v)^2$.
\item[(iv)] the triple $\big(I,T(z)R(x,y),T(z)L(u,v)\big)$ is an isotopism
from $G$ to $G'$.
\end{description}
\item[(b)] if $P(z,R(x,y)):=0$ then,
\begin{description}
\item[(i)] $G$ is an A$_\rho$-loop and $|R(x,y)|=2$.
\item[(ii)] $\textrm{Inn}_\lambda (G)=\big<R(x,y)L_e'~:~x,y\in \Omega\big>$.
\item[(iii)] the triple $(R(x,y)T(z),I,R(x,y)L(u,v))$ is an isotopism
from $G$ to $G'$.
\end{description}
\item[(c)] if $P(z,R(x,y)):=0$ then,
\begin{description}
\item[(i)] $G$ is an A$_\lambda$-loop and $|R(x,y)|=2$.
\item[(ii)] $\textrm{Inn}_\lambda (G)=\big<R(x,y)L_e'~:~x,y\in \Omega\big>$.
\item[(iii)] the triple $\big(L(x,y)T(z),L(x,y)R(u,v),I\big)$ is an isotopism
from $G$ to $G'$.
\end{description}
\end{description}
\end{enumerate}
\end{myth}
\textit{Proof}\\
This is proved using Theorem~\ref{main:result1} and
Corollary~\ref{main:result2}.

\begin{mycor}\label{2:2}
Let $G=(\Omega ,\cdot )$ and $G'=(\Omega ,\circ )$ be two distinct
isotopic loops.
\begin{enumerate}
\item Under the triple $\big(R(x,y),L(u,v),T(z)\big)$, if $P(z,R(x,y)):=0$, $P(z,L(x,y)):=0$ and $P(z,T(x)):=0$
then,
\begin{description}
\item[(a)] $G$ is an A-loop and $|R(x,y)|=|L(x,y)|=|T(x)|=2$.
\item[(b)] $T(z)L(x,y)=L(x,y)T(z)$.
\item[(c)] $\textrm{Inn}_\mu (G)=\big<L(x,y)L_e'~:~x,y\in \Omega\big>$ and $\textrm{Inn}_\lambda (G)=\big<T(x)L_e'~:~x\in
\Omega\big>$.
\item[(d)] the triples $\big(I,R(x,y)L(u,v),R(x,y)T(z)\big)$,$\big(L(u,v)R(x,y),I,L(u,v)T(z)\big)$ and
$\big(T(z)R(x,y),T(z)L(u,v),I\big)$ are isotopisms from $G$ to $G'$.
\end{description}
\item Under the triple $\big(L(x,y),R(u,v),T(z)\big)$, if $P(z,R(x,y)):=0$, $P(z,L(x,y)):=0$ and $P(z,T(x)):=0$
then,
\begin{description}
\item[(a)] $G$ is an A-loop and $|R(x,y)|=|L(x,y)|=|T(x)|=2$.
\item[(b)] $T(z)R(x,y)=R(x,y)T(z)$.
\item[(c)] $\textrm{Inn}_\mu (G)=\big<R(x,y)L_e'~:~x,y\in \Omega\big>$ and $\textrm{Inn}_\rho (G)=\big<T(x)L_e'~:~x\in
\Omega\big>$.
\item[(d)] the triples $\big(I,L(x,y)R(u,v),L(x,y)T(z)\big)$, $\big(R(x,y)L(u,v),I,R(x,y)T(z)\big)$ and
$\big(T(z)L(x,y),T(z)R(u,v),I\big)$ are
isotopisms from $G$ to $G'$.
\end{description}
\item Under the triple $\big(T(z),R(x,y),L(u,v)\big)$, if $P(z,R(x,y)):=0$, $P(z,L(x,y)):=0$ and $P(z,T(x)):=0$
then,
\begin{description}
\item[(a)] $G$ is an A-loop and $|R(x,y)|=|L(x,y)|=|T(x)|=2$.
\item[(b)] $R(x,y)L(u,v)=L(u,v)R(x,y)$.
\item[(c)] $\textrm{Inn}_\lambda (G)=\big<R(x,y)L_e'~:~x,y\in \Omega\big>$ and $\textrm{Inn}_\rho (G)=\big<L(x,y)L_e'~:~x\in
\Omega\big>$.
\item[(d)] the triples $\big(I,T(z)R(x,y),T(z)L(u,v)\big)$, $(R(x,y)T(z),I,R(x,y)L(u,v))$ and
$\big(L(x,y)T(z),L(x,y)R(u,v),I\big)$ are isotopisms from $G$ to
$G'$.
\end{description}
\end{enumerate}
\end{mycor}
\textit{Proof}\\
This follows directly from Theorem~\ref{2:1}.


\begin{thebibliography}{99}
\bibitem{phd41} R. H. Bruck, {\it A survey of binary systems}. Springer-Verlag, Berlin-G\"ottingen-Heidelberg, 1966.
\bibitem{phd40} R. H. Bruck and L. J. Paige, {\it Loops whose
inner mappings are automorphisms}. The annuals of Mathematics,
\textbf{63}(1956), 2, 308--323.
\bibitem{phd39} O. Chein, H. O. Pflugfelder and J. D. H. Smith, {\it Quasigroups and loops : Theory and applications}. Heldermann Verlag, 1990.
\bibitem{phd108} P. Cs\"org\H o, {\it Extending the
structural homomorphism of LCC loops}. Comment. Math. Univ.
Carolinae \textbf{46}(2005), 3, 385--389.
\bibitem{phd106} P. Cs\"org\H o and  A. Dr\'apal, {\it Left
conjugacy closed loops of nilpotency class 2}. Results Math.
\textbf{47}(2005), 3-4, 242--265.
\bibitem{phd103} P. Cs\"org\H o and T. Kepka , {\it On loops
whose inner permutations commute}. Comment. Math. Univ. Carolinae
\textbf{45}(2004), 2, 213--221.
\bibitem{phd164} P. Cs\"org\H o and M. Niemenmaa, {\it Solvability
conditions for loops and groups}. J. Alg. \textbf{232}(2000), 1,
336--342.
\bibitem{phd49} J. Dene and A. D. Keedwell, {\it Latin squares and their applications}. The English University press Lts, 1974.
\bibitem{phd14} A. Dr\'apal, {\it A-loops close to code loops are groups}. Comment. Math. Univ. Carolinae \textbf{41}(2000), 2, 245--249.
\bibitem{phd37} A. Dr\'apal, {\it Conjugacy closed loops and
their multiplication groups}. J. Alg. \textbf{272}(2004), 1,
838--850.
\bibitem{phd38} A. Dr\'apal, {\it Structural interactions of conjugacy closed loops}. Trans. Amer. Math.
Soc. To appear.
\bibitem{phd98} A. Dr\'apal, {\it On multiplicatipon groups of left conjugacy closed
loops}. Comment. Math. Univ. Carolinae \textbf{45}(2004), 223--236.
\bibitem{phd139} A. Dr\'apal, {\it On extraspecial left conjugacy closed loops}. J. Alg. \textbf{302}(2006), 2, 771--792.
\bibitem{phd107} A. Dr\'apal, {\it On left conjugacy closed loops with a nucleus of index two}. Abh. Math. Sem. Univ. Hamburg. \textbf{74}(2004), 205--221.
\bibitem{phd99} A. Dr\'apal, {\it Multiplication groups of finite free loops that fix at most two points}. J. Alg. \textbf{235}(2001), 1, 154--175.
\bibitem{phd100} A. Dr\'apal, {\it Multiplication groups of free loops
I}. Czechoslovak. Math. J. \textbf{46}(1996), 121--131.
\bibitem{phd101} A. Dr\'apal, {\it Multiplication groups of free loops
II}. Czechoslovak. Math. J. \textbf{46} (1996), 201--220.
\bibitem{phd102} A. Dr\'apal, {\it Orbits of inner mapping
groups}. Monatsh. Math. \textbf{134}(2002), 191--206.
\bibitem{phd104} A. Dr\'apal and P. Cs\"org\H o
, {\it On left conjugacy closed loops in which the left
multiplication group is normal}. Manuscript.
\bibitem{phd129} A. Dr\'apal and T. Kepka, {\it
Multiplication groups of quasigroups and loops I}. Acta Univ. Carol.
Math. Phys. \textbf{34}(1993), 1, 85--99.
\bibitem{phd130} A. Dr\'apal, T. Kepka and O. Mar\v s\'alek, {\it
Multiplication groups of quasigroups and loops II}. Acta Univ.
Carol. Math. Phys. \textbf{35}(1994), 1, 9--29.
\bibitem{phd91} E. G. Goodaire and D. A. Robinson, {\it A
class of loops which are isomorphic to all loop isotopes}. Can. J.
Math. \textbf{34}(1982), 662--672.
\bibitem{phd48} E. G. Goodaire and D. A. Robinson, {\it Some special conjugacy closed loops}. Canad. Math. Bull. \textbf{33}(1990), 73--78.
\bibitem{phd42} E. G. Goodaire, E. Jespers and C. P. Milies, {\it Alternative loop rings}. NHMS(184), Elsevier, 1996.
\bibitem{phd131} T. Kepka, {\it Multiplication groups of some
quasigroups}. Colloq. Math. Soc. J. Bolyai \textbf{29}(1977),
459--465.
\bibitem{phd132} T. Kepka, {\it On the abelian inner
permutation groups of loops}. Comm. Alg. \textbf{26}(1998),
857--861.
\bibitem{phd133} T. Kepka and M. Niemenmaa, {\it On loops
with cyclic inner mapping groups}. Arch. Math. \textbf{60}(1993),
233--236.
\bibitem{phd36} M. K. Kinyon, K. Kunen, {\it The structure of
extra loops}. Quasigroups and Related Systems \textbf{12}(2004),
39--60.
\bibitem{phd35} M. K. Kinyon, K. Kunen, J. D. Phillips, {\it
Diassociativity in conjugacy closed loops}. Comm. Alg.
\textbf{32}(2004), 767--786.
\bibitem{phd47} M. K. Kinyon, K. Kunen, {\it
Power-associative conjugacy closed loops}. J. Alg.
\textbf{304}(2006), 2, 679--711.
\bibitem{phd13} M. K. Kinyon, K. Kunen, J. D. Phillips, {\it Every diassociative A-Loop is Moufang}. Proc. Amer. Math. Soc. \textbf{130}(2001), 3, 619--624.
\bibitem{phd78} K. Kunen, {\it The structure of conjugacy closed loops}. Trans. Amer. Math. Soc. \textbf{352}(2000), 2889--2911.
\bibitem{phd135} M. Niemenmaa, {\it On the structure of the inner
mapping groups of loops}. Comm. Alg. \textbf{24}(1996), 135--142.
\bibitem{phd136} M. Niemenmaa, {\it On finite loops whose inner
mapping groups are abelian}. Bull. Aust. Math. Soc.
\textbf{65}(2002), 477--484.
\bibitem{phd162} M. Niemenmaa, {\it On loops which have
dihedral $2$-groups as inner mappings groups}. Bull. Aust. Math.
Soc. \textbf{52}(1995), 153--169.
\bibitem{phd134} M. Niemenmaa and T. Kepka, {\it
On multiplication groups of loops}. J. Alg. \textbf{135}(1990),
112--122.
\bibitem{phd110} J. M. Osborn, {\it A theorem on A-loops}. Proc. Amer. Math. Soc. \textbf{9}(1959), 347--349.
\bibitem{phd3} H. O. Pflugfelder, {\it Quasigroups and loops : Introduction}. Sigma series in Pure Math. 7, Heldermann Verlag, Berlin, 1990.
\bibitem{phd123} J. D. Phillips, {\it On Moufang A-loops}. Comment. Math. Univ. Carolinae \textbf{41}(2000), 2, 371--375.
\bibitem{phd151} J. D. Phillips, {\it A short basis for the variety of WIP
PACC-loops}. Quasigroups and Related Systems \textbf{14}(2006), 1,
73--80.
\bibitem{phd75} W. B. Vasantha Kandasamy, {\it Smarandache
loops}. Department of Mathematics, Indian Institute of Technology,
Madras, India, 2002.
\bibitem{phd163} A. Vesanen, {\it Solvable groups and loops},
J. Alg. \textbf{180}(1996), 862--876.
\end{thebibliography}
\end{document}